\documentclass[11pt]{amsart}
\usepackage[utf8]{inputenc}
\usepackage[T1]{fontenc}
\usepackage{lmodern}
\usepackage[margin=1in]{geometry}
\usepackage{amsmath,amssymb,amsthm,colonequals,mathrsfs,mathtools}
\usepackage{array,enumitem,booktabs,longtable,xcolor}
\usepackage{microtype}
\usepackage[colorlinks,linkcolor=blue,citecolor=blue,urlcolor=blue,bookmarksdepth=2]{hyperref}
\allowdisplaybreaks
\AtBeginDocument{\def\MR#1{}}
\theoremstyle{plain}
\newtheorem{thm}{Theorem}[section]
\newtheorem{lem}[thm]{Lemma}
\newtheorem{cor}[thm]{Corollary}
\newtheorem{prop}[thm]{Proposition}

\theoremstyle{definition}

\newtheorem{ex}[thm]{Example}
\numberwithin{equation}{section}

\def\Q{{\mathbb Q}}
\def\R{{\mathbb R}}
\def\Z{{\mathbb Z}}
\def\C{{\mathbb C}}

\def\Aut{\mathop{\mathrm{Aut}}\nolimits}
\def\Gal{\mathop{\mathrm{Gal}}\nolimits}
\def\Lie{\mathop{\mathrm{Lie}}\nolimits}
\def\Hom{\mathop{\mathrm{Hom}}\nolimits}

\def\GL{\mathop{\mathrm{GL}}\nolimits}
\def\SL{\mathop{\mathrm{SL}}\nolimits}

\def\det{\mathop{\mathrm{det}}\nolimits}
\def\dim{\mathop{\mathrm{dim}}\nolimits}

\def\Ind{\mathop{\mathrm{Ind}}\nolimits}

\def\Hom{\mathop{\mathrm{Hom}}\nolimits}
\def\diag{\mathop{\mathrm{diag}}\nolimits}

\def\Spin{\mathop{\mathrm{Spin}}\nolimits}
\def\exp{\mathop{\mathrm{exp}}\nolimits}
\def\Sp{\mathop{\mathrm{Sp}}\nolimits}
\def\SO{\mathop{\mathrm{SO}}\nolimits}

\def\Sym{\mathop{\mathrm{Sym}}\nolimits}

\def\U{\mathrm{U}}
\def\O{\mathrm{O}}

\newcommand{\bbA}{\mathbb{A}}

\newcommand{\bbC}{\mathbb{C}}

\newcommand{\bbF}{\mathbb{F}}

\newcommand{\bbG}{\mathbb{G}}

\newcommand{\bbQ}{\mathbb{Q}}

\newcommand{\bbR}{\mathbb{R}}

\newcommand{\bbZ}{\mathbb{Z}}

\newcommand{\calD}{\mathcal{D}}

\newcommand{\calJ}{\mathcal{J}}

\newcommand{\calO}{\mathcal{O}}

\newcommand{\calS}{\mathcal{S}}

\newcommand{\calX}{\mathcal{X}}

\newcommand{\frakf}{\mathfrak{f}}

\newcommand{\frakH}{\mathfrak{H}}

\newcommand{\frakp}{\mathfrak{p}}

\newcommand{\bfi}{\mathbf{i}}

\newcommand{\Mat}{\mathrm{Mat}}

\newcommand{\sgn}{\mathrm{sgn}}

\newcommand{\fini}{\mathrm{fin}}
\newcommand{\bs}{\backslash}

\newcommand{\vep}{\varepsilon}
\newcommand{\isom}{\cong}

\newcommand{\Irr}{\mathrm{Irr}}

\newcommand{\Ad}{\mathrm{Ad}}

\newcommand{\BC}{\mathrm{BC}}

\newcommand{\As}{\mathrm{As}}

\def\val{\mathop{\mathrm{val}}\nolimits}

\allowdisplaybreaks[4]

\newcommand{\defeq}{\vcentcolon=}

\newcommand{\Cl}{\mathrm{Cl}}

\newcommand{\DD}{\mathcal D}

\newcommand{\GG}{\Gamma}
\newcommand{\sn}{\operatorname{sn}}
\newcommand{\age}{\operatorname{age}}

\begin{document}
\title[Kodaira dimension of modular varieties of tube type]{The Kodaira dimension of modular varieties of tube type}
\author{Shuji Horinaga}
\address{NTT Institute for Fundamental Mathematics, NTT Inc., Japan}
\email{syuuji.horinaga@ntt.com, shorinaga@gmail.com}
\author{Yota Maeda}
\address{Fachbereich Mathematik, Technische Universit\"at Darmstadt, Germany / Mathematical Institute, Tohoku University, Japan}
\email{y.maeda.math@gmail.com}
\date{}

\begin{abstract}
We study the Kodaira dimension of arithmetic quotients of classical domains of tube type.
For the unitary groups of split Hermitian spaces and for the groups over definite quaternion algebras, we give bounds for general type that hold for every arithmetic subgroup.
We also prove that every arithmetic quotient of the Siegel upper half space of genus $12$
is of general type, which extends the result by Tai.
In the course of the proof, we construct low weight cusp forms for every arithmetic subgroup whose partial $L$-functions coincide with those of the Ikeda and Miyawaki lifts by means of automorphic representation theory.
For domains of type DIII, including those that are not of
tube type, we prove a sharp bound for canonicity of singularities on toroidal compactifications. Our proof method also extends a theorem of Ma on the
Kodaira dimension of orthogonal modular varieties to rational projective similitudes.
\end{abstract}
\maketitle

\section{Introduction}

The Kodaira dimension of modular varieties has long been studied in connection
with moduli spaces of polarized abelian varieties and of K3 surfaces.
For Siegel modular varieties we refer to Tai, Freitag, and Mumford
\cite{tai1982kodaira,freitag1983Siegelsche,Mumford1983Siegel,Tai1995Polarizations},
and for orthogonal modular varieties to Gritsenko, Hulek, and Sankaran
\cite{Gritsenko2007kodaira} and Ma \cite{ma2018kodaira}.
In both settings, a cusp form of sufficiently small weight can be used
to construct pluricanonical forms on a smooth projective model.

The classical irreducible domains of tube type are associated with
$\mathrm{Sp}_{2n}(\mathbb R)$, $\mathrm U(n,n)$, $\mathrm{SO}^*(4n)$,
and $\mathrm{SO}_0(2,n)$ \cite{1994_FK_analysis_symm_cone}.
We study the unitary and quaternionic families by this method and obtain
further results for the symplectic and orthogonal families.
One motivation comes from a question of van der Geer about the geometric
meaning of the Ikeda lift in genus $12$ of the discriminant form
$\Delta$ of weight $12$ for $\mathrm{SL}_2(\mathbb Z)$
\cite{vanDerGeer2008Siegel}.

\subsection{Kodaira dimension}
\label{subsec:intro-kodaira}

Let $A$ be $\Q$, an imaginary quadratic field $E$ of discriminant $-D$,
or a definite quaternion algebra $B$ over $\Q$. We equip $A$ with the
identity, complex conjugation, or the main involution, respectively.
For a matrix $x$, let $x^*$ be its transpose with this involution applied
to each entry. On the right $A$-module $A^{2n}$
we take the form $(x,y)\mapsto x^*J_ny$, where
$J_n=\left(\begin{smallmatrix}0&I_n\\-I_n&0\end{smallmatrix}\right)$.
Let $H_{A,n}$ be the identity component of the group of isometries over $\Q$.
Write $\frakH_{A,n}$ for its symmetric domain and $k_{\mathrm{can}}$ for
the weight of the canonical bundle.
\[
\begin{array}{c|c|c|c}
 A & H_{A,n}(\R) & \frakH_{A,n} & k_{\mathrm{can}}\\ \hline
 \Q & \Sp_{2n}(\R) & \Sp_{2n}(\R)/\mathrm U(n) & n+1\\
 E & \mathrm U(n,n) & \mathrm U(n,n)/(\mathrm U(n)\times\mathrm U(n)) & 2n\\
 B & \mathrm{SO}^*(4n) & \mathrm{SO}^*(4n)/\mathrm U(2n) & 2n-1
\end{array}
\]
Here $\frakH_{\Q,n}$ is the Siegel upper half space of genus $n$.
For $A=B$ the group $H_{B,n}$ consists of the isometries of reduced norm one.
For an arithmetic subgroup $\Gamma\subset H_{A,n}(\Q)$, put
$X_\Gamma=\Gamma\backslash\frakH_{A,n}$. We write $\kappa(X_\Gamma)$
and $p_g(X_\Gamma)$ for the Kodaira dimension and geometric genus
of a smooth projective model of $X_\Gamma$. For $A=E$, let $c_E$ be $8$ if
$D=3$, $6$ if $D\in\{4,7\}$, $4$ if $D\in\{8,11,15,20, 23\}$, and $2$ otherwise.

\begin{thm}[Corollary~\ref{cor:kodaira-degrees}]
\label{mainthm:kodaira}
For every arithmetic subgroup $\Gamma\subset H_{A,n}(\Q)$ the following hold.
\begin{enumerate}
\item If $A=\Q$, then $X_\Gamma$ is of general type for $n=12$.
For $n\in\{11,15,17\}$ it satisfies $p_g(X_\Gamma)\geq1$ and $\kappa(X_\Gamma)\geq0$.
\item If $A=E$, then $X_\Gamma$ is of general type for odd $n\geq13$
and for even $n>c_E$. For $n=11$ or $n=c_E$ it satisfies
$p_g(X_\Gamma)\geq1$ and $\kappa(X_\Gamma)\geq0$.
\item If $A=B$, then $X_\Gamma$ is of general type for $n\geq7$.
For $n=6$ it satisfies
$p_g(X_\Gamma)\geq1$ and $\kappa(X_\Gamma)\geq0$.
\end{enumerate}
\end{thm}

For the moduli space $\mathcal A_g$ of principally polarized abelian
varieties, Tai \cite{tai1982kodaira} proved general type for $g\geq9$,
Freitag \cite{freitag1983Siegelsche} for $g\geq8$, and Mumford
\cite{Mumford1983Siegel} for $g\geq7$.
Wang \cite{Wang2026A6} proved that $\kappa(\mathcal A_6)=0$.
For each fixed type of polarization, Tai \cite{Tai1995Polarizations}
proved general type for $g\geq16$; see also \cite{HulekSankaranSurvey}.
For abelian surfaces, more explicit results on the birational geometry of
these moduli spaces are available for specific polarization types and level
structures \cite{OGrady1989Abelian,HulekSankaran1994Abelian,GritsenkoSankaran1996Polarisation,Sankaran1997Polarised,Hulek2001Igusa,GrossPopescu2001Abelian,Erdenberger2004Abelian}.
Part~(1) gives general type for every type of polarization in genus $12$
and a nonzero differential of top degree for every arithmetic subgroup in genera $11$, $15$, and $17$.

In the unitary case, Horinaga, Maeda, and Yamauchi
\cite{2025_HMY_Kodaira_dimension_ball_quotients} studied ball quotients of signature $(1,n)$. For signature $(p,q)$, the
criterion for general type in \cite[Theorem~2.3]{2026_Horinaga-Maeda_sing_unitary}
assumes $p+q\equiv1\pmod{|\mathcal O_E^\times|}$. Since $|\mathcal O_E^\times|$ is even, this condition excludes the signature $(n,n)$ of part~(2).

The Kodaira dimension of quaternionic modular varieties of higher rank, as in
part~(3), is less well understood. For comparison, Voight
\cite{Voight2009Shimura} classified the Shimura curves of genus at most two
arising from the units of reduced norm one in Eichler orders, and thus
determined the exceptions to general type in this family of curves.

The extension of pluricanonical forms to a
resolution requires a bound for the singularities. We prove the following
theorem for the whole family of domains of type DIII.
Let $\DD^{\mathrm{DIII}}_m=\mathrm{SO}^*(2m)/\mathrm U(m)$.
Its complex dimension is $m(m-1)/2$, and it is of tube type exactly
when $m$ is even. In particular, $\frakH_{B,n}\simeq\DD^{\mathrm{DIII}}_{2n}$.

\begin{thm}[Theorem~\ref{thm:DIII-compactification}]\label{thm:DIII}
Let $m\geq7$ and fix an absolutely simple adjoint group $\bbG$ over $\Q$ with $\bbG(\bbR)^\circ = \Aut(\DD^{\mathrm{DIII}}_m)^\circ$.
For every arithmetic subgroup $\GG$ of $\bbG(\bbQ) \cap \bbG(\bbR)^\circ$, the quotient $\GG\backslash\DD^{\mathrm{DIII}}_m$ admits
a projective toroidal compactification with canonical singularities.
\end{thm}

Example~\ref{ex:sharpness} shows that the bound is sharp as the rational
form and the arithmetic subgroup vary.
For domains of type DIII we estimate the age of an element of finite order
on the exterior square of the standard representation. We then compare its
action in the interior and at the boundary of a toroidal compactification.
Let $V$ be a rational quadratic space of signature $(2,n)$ and let $L$ be
an integral lattice in $V$. Ma \cite[Corollary~1.8(1), (2)]{ma2018kodaira}
proved that the orthogonal modular variety for $\mathrm{SO}^+(L)$ is of
general type for even $n\geq22$ and has nonnegative Kodaira dimension for
$n=16,20$. In Theorem~\ref{thm:orthogonal}, we use the same comparison and
descent through a central isogeny to extend these bounds to arithmetic
subgroups of $\mathrm{PGSO}(V)(\Q)$ preserving the chosen component of the
domain.

\subsection{Construction of cusp forms}

\begin{thm}\label{thm:intro-cusp-forms}
Let $(n,\ell)$ be a degree and weight listed for $A$ in Table~\ref{tab:cusp-weights} below.
For every arithmetic subgroup $\Gamma\subset H_{A,n}(\Q)$ there is a
nonzero cusp form of weight $\ell$ for $\Gamma$ on $\frakH_{A,n}$.
\end{thm}

\begin{table}[!ht]
\caption{Degrees and weights of the cusp forms}\label{tab:cusp-weights}
\centering
\(
\begin{array}{c|l|c}
 A & \text{degree }n & \text{weight }\ell\\ \hline
 \Q & n\equiv0\pmod4,\quad n\geq12 & n/2+6\\
    & n\equiv2\pmod4,\quad n\geq18 & n/2+9\\
    & n\equiv3\pmod4,\quad n=11\ \text{or}\ n\geq19 & (n+13)/2\\
    & n=15 & 16\\
    & n\equiv1\pmod4,\quad n\geq17 & (n+19)/2\\ \hline
 E & n\text{ odd},\quad n\geq11 & n+11\\
   & n\text{ even},\quad n\geq c_E & n+c_E\\ \hline
 B & n\geq3 & n+5
\end{array}
\)
\end{table}

We use the theory of automorphic representations to pass
from lifts of elliptic modular forms to cusp forms for the prescribed
group $\Gamma$. For $A=\Q$ and $n$ even we use the Ikeda lift \cite{ikeda2001lifting}.
For $n$ odd, we use Arthur's multiplicity formula for the parameter of the Miyawaki lift, without assuming that the integral defining the lift
is nonzero \cite{2006_Ikeda_Miyawaki-lift}. 
For $A=E$ we apply Arthur's multiplicity
formula \cite{2015_Mok} to the parameters of the lifts studied by Ikeda \cite{2008_Ikeda}. 
In both families we choose local representations
with vectors fixed by the prescribed compact open subgroups and exclude the residual spectrum by investigating the lowest $K_\infty$-type at infinity applying \cite{2008_Moeg_noncusp}. 
For $A=B$ we use Yamana's description of the local components of his lift, including those at the
ramified primes \cite[Theorem~4.2]{Yamana}. In the last two families
we prove that the restrictions to all connected components have the
same Petersson norm. This gives a cusp form on the prescribed
component without a hypothesis on the class number.

\subsection*{Acknowledgments}
The authors thank Klaus Hulek for his insightful comments.
The first author is partially supported by JSPS KAKENHI Grant Number 23K12965.
The second author is partially supported by the Alexander von Humboldt Foundation through a Humboldt research fellowship, by Deutsche Forschungsgemeinschaft (DFG, German Research Foundation) through the Collaborative Research Centre TRR 326 Geometry and Arithmetic of Uniformized Structures, project number 444845124, and by Waseda University Grant for Special Research Projects (Project number: 2026R-022).
The authors used AI in the preparation of the manuscript. The authors are responsible for the mathematical content and the final manuscript.

\section*{Notation and conventions}
For $A$ as in Section~\ref{subsec:intro-kodaira}, put
$A_v=A\otimes_\Q\Q_v$ for each place $v$ of $\Q$. For $A=E$ we fix an
embedding $E\hookrightarrow\C$. Multiplying $J_n$ by a nonzero purely
imaginary element of $E$ gives a split Hermitian form.
A Hermitian symmetric domain is of tube type if it is biholomorphic to $V+\sqrt{-1}\Omega$ for an open convex cone $\Omega$ in a real vector space $V$. 
For $\frakH_{A,n}$ we take $V=\{x\in\Mat_n(A_\infty)\mid x^*=x\}$ and let $\Omega$ be the cone of positive definite elements of $V$. We write $K_\infty$ for the stabilizer of $\bfi=\sqrt{-1}\,I_n$ in $H_{A,n}(\R)$.
An element $g=\left(\begin{smallmatrix}a&b\\c&d\end{smallmatrix}\right)$
of $H_{A,n}(\R)$ acts by $gz=(az+b)(cz+d)^{-1}$.
Let $\nu$ be the reduced norm over the center of $A$, which is the determinant
for $A=\Q,E$. Following \cite[Section~1]{Yamana}, we set
$\calJ_k(g,z)=\nu(cz+d)^k$ for $A=\Q,B$. For $A=E$ we take $k$ even
and set $\calJ_k(g,z)=\det(g)^{-k/2}\nu(cz+d)^k$.
For $A=B$ our weight $k$ corresponds to weight $2k$ in
\cite[p.~2485 and Definition~1.1]{Yamana}.
A modular form of weight $k$ for $\Gamma$ is a holomorphic function
$F$ on $\frakH_{A,n}$ that is holomorphic at the rational cusps and
satisfies $F(\gamma z)=\calJ_k(\gamma,z)F(z)$ for every $\gamma\in\Gamma$.
We write $M_k(\Gamma)$ for this space and $S_k(\Gamma)$ for the subspace
of cusp forms.
Let $\omega$ be a nonzero constant differential of top degree in tube
coordinates. Then $g^*\omega=\calJ_{k_{\mathrm{can}}}(g,z)^{-1}\omega$.
For $A=B$ the Jacobian is $\nu(cz+d)^{-(2n-1)}$
\cite[Chapter~II, Theorem~1.10, p.~52]{Krieg}
(see also \cite[Lemmas~5.2 and~5.3]{FreitagSalvatiManni}).

\section{Automorphic forms and connected components}

Write $\bbA$ and $\bbA_{\mathrm{fin}}$ for the adeles and finite adeles of $\Q$.
Fix a Haar measure on $H_{A,n}(\bbA)$ and write $dg$ for the quotient measure.
For square integrable automorphic forms we use the Petersson inner product
\[
 \langle\varphi_1,\varphi_2\rangle
 =\int_{H_{A,n}(\Q)\backslash H_{A,n}(\bbA)}
   \varphi_1(g)\overline{\varphi_2(g)}\,dg,
 \qquad \lVert\varphi\rVert^2=\langle\varphi,\varphi\rangle.
\]

\subsection{Fixed vectors in induced representations}
Let $P=MN$ be a parabolic subgroup of a connected reductive group $L$ over a
nonarchimedean local field $F$. Write $\Ind_{P(F)}^{L(F)}$ for normalized
parabolic induction and $\delta_P$ for the modulus character of $P(F)$.
We shall apply the following lemma to characters of the Levi that factor through the valuation of a determinant or a reduced norm.

\begin{lem}\label{lem:fixed}
With the above notation, suppose that a character $\xi$ of $M(F)$ is trivial on every compact subgroup.
For each compact open subgroup $K\subset L(F)$, the space $\bigl(\Ind_{P(F)}^{L(F)}\xi\bigr)^K$ is nonzero.
\end{lem}
\begin{proof}
The image of $P(F)\cap K$ in $M(F)$ is compact. Both $\xi$ and $\delta_P^{1/2}$ are trivial on this image.
Let $f$ be a function whose support is $P(F)K$ such that $f(mnk)=\delta_P(m)^{1/2}\xi(m)$ for $m\in M(F)$, $n\in N(F)$ and $k\in K$.
It is a well-defined nonzero $K$-fixed vector in the induced representation.
\end{proof}

For $n=\sum_i n_i$ and a representation $\tau_i$ of $\GL_{n_i}(\Q_p)$, let $\tau_1\times\cdots\times\tau_r$ be the normalized induction of $\GL_n(\Q_p)$ from the standard parabolic with Levi subgroup $\prod_i\GL_{n_i}$. 
The same notation applies over a finite extension of $\Q_p$.
For a rational parabolic subgroup $Q$ of $H_{A,n}$ with blocks of sizes
$n_1,\ldots,n_r$, put $n_0=n-\sum_{i=1}^r n_i$ and take $H_{A,0}=1$.
The Levi subgroup of $Q(\Q_v)$ is isomorphic to
$\prod_{i=1}^{r}\GL_{n_i}(A_v)\times H_{A,n_0}(\Q_v)$.
Take representations $\tau_i$ of $\GL_{n_i}(A_v)$ and $\pi_0$ of $H_{A,n_0}(\Q_v)$.
Define
\[
\tau_1 \times \cdots \times \tau_r \rtimes \pi_0 \defeq \Ind_{Q(\bbQ_v)}^{H_{A,n}(\bbQ_v)} \tau_1 \boxtimes \cdots \boxtimes \tau_r \boxtimes \pi_0.
\]

\subsection{Nonvanishing on each connected component}
\label{subsec:classical-symplectic-unitary}

Let $K\subset H_{A,n}(\bbA_\fini)$ be a compact open subgroup.
Write $\Cl_K$ for the group indexing the connected components of $H_{A,n}(\Q)\backslash H_{A,n}(\bbA)/(K_\infty K)$. 
It is trivial for $A=\Q$ by strong approximation. 
For $H_{E,n}$ put $T=\ker(N_{E/\Q}:\operatorname{Res}_{E/\Q}\mathbb G_m\to\mathbb G_m)$.
The determinant and strong approximation for the special unitary group give $\Cl_K=T(\Q)\backslash T(\bbA_\fini)/\det K$.
For $A=B$ and $n\geq3$, let $\sn$ be the spinor norm defined by the connecting map for the spin covering $\operatorname{Spin}(H_{B,n})\to H_{B,n}$. 
At each finite prime $p$, the spinor norm $\sn\colon H_{B,n}(\Q_p)\to\Q_p^\times/\Q_p^{\times2}$ is surjective by \cite[Proposition 6.5]{2011_Yamana_degenerate_PS}. Its image at the real place is trivial \cite[Remark~6.6]{2011_Yamana_degenerate_PS}.
By exactness for the spin covering and the Hasse principle for the spin group \cite[Theorem~5.12.31]{PoonenRationalPoints}, a class in $\bbQ^\times/ \bbQ^{\times2}$ is a spinor norm precisely when it is a spinor norm at every place. 
To identify the double cosets, take two representatives in $H_{B,n}(\bbA_\fini)$ whose spinor norms agree modulo $\sn(H_{B,n}(\Q))\sn(K)$. After adjustment by a rational point and an element of $K$, their quotient lifts to $\Spin(H_{B,n})(\bbA_\fini)$. The spin group is simply connected, absolutely almost simple, and noncompact at infinity.
Therefore $\Spin(H_{B,n})(\bbQ)$ is dense in $\Spin(H_{B,n})(\bbA_\fini)$. 
We obtain
\[
 \sn(H_{B,n}(\Q))=\Q_{>0}^{\times}/\Q^{\times2},\qquad
 \Cl_K\simeq\Q_{>0}^{\times}\backslash\bbA_\fini^\times/\widetilde{\sn(K)} \isom \bbQ^\times \bs \bbA^\times/\bbA^{\times2}\bbR_{>0}^\times\widetilde{\sn(K)},
\]
where $\widetilde{\sn(K)}$ is the inverse image of $\sn(K)$ in $\bbA_\fini^\times$.

For $x \in \Cl_K$, fix a representative $g_x \in H_{A,n}(\bbA_\fini)$ such that $g_1 = 1$.
We have the decomposition
\[
H_{A,n}(\bbQ) \bs H_{A,n}(\bbA)/K_\infty K 
= \bigsqcup_{x \in \Cl_K} \Gamma_x \bs H_{A,n}(\bbR)/K_\infty
= \bigsqcup_{x \in \Cl_K} \Gamma_x \bs \frakH_{A,n}, 
\]
where $\Gamma_x = H_{A,n}(\bbQ) \cap g_x K g_x^{-1}$.
Choose a compact open subgroup
$K_{\Gamma,p}$ containing $\Gamma$ such that $K_{\Gamma,p}$ is hyperspecial for almost all $p$. 
Write $K_{\Gamma} = \prod_p K_{\Gamma,p}$.
At this level $\Gamma\subset\Gamma_1$,
so every cusp form for $\Gamma_1$ is a cusp form for $\Gamma$.
For an automorphic form $\varphi$ let $\varphi|_x$ be its restriction to
$H_{A,n}(\Q)g_x H_{A,n}(\R)K$.
For a character $\chi$ of $\Cl_K$, let $\widetilde\chi$ be its pullback to $H_{A,n}(\bbA)$
by the map to $\Cl_K$, given by the determinant for $H_{E,n}$ and by the spinor norm for $H_{B,n}$.

\begin{lem}\label{prop:common-components}
Let $K\subset H_{A,n}(\bbA_\fini)$ be compact and open, and let $\pi$ be a discrete
automorphic representation with $\pi^K\ne0$. Suppose that
$\pi\not\simeq\pi\otimes\widetilde\chi$ for every nontrivial character $\chi$ of $\Cl_K$.
Then $\lVert\varphi|_x\rVert^2=\lVert\varphi\rVert^2/|\Cl_K|$ for every
nonzero $\varphi\in\pi^K$ and every $x\in\Cl_K$.
\end{lem}
\begin{proof}
    Let $\chi$ be a nontrivial character of $\Cl_K$ and $\varphi_\chi(g) = \widetilde\chi(g) \varphi(g)$.
    The automorphic form $\varphi_\chi$ generates $\pi \otimes \widetilde{\chi}$.
    Since $\pi$ and $\pi\otimes\widetilde{\chi}$ are inequivalent, we have
    \[
    0 = \langle \varphi_\chi, \varphi\rangle = \sum_{x \in \Cl_K} \chi(x)\lVert\varphi|_x\rVert^2.
    \]
Hence, Fourier transform on $\Cl_K$ for the function $x \mapsto \lVert \varphi|_x\rVert^2$ shows that $\lVert\varphi|_x\rVert^2$ is independent of $x$.
\end{proof}

Let $\frakp_-$ be the subspace of the complexified Lie algebra of $H_{A,n}(\R)$
corresponding to the antiholomorphic tangent space at $\bfi$.
Define the character $\tau_\ell(k)=\calJ_\ell(k,\bfi)^{-1}$ of $K_\infty$.
Let $\pi_\infty(\ell)$ be the lowest weight representation with lowest $K_\infty$-type $\tau_\ell$.
The space of automorphic forms fixed by $K$, transforming under $K_\infty$
by $\tau_\ell$, and annihilated by $\frakp_-$ is naturally isomorphic to
$\bigoplus_{x\in\Cl_K}M_\ell(\Gamma_x)$. The $x$-component of the image of
$\varphi$ takes the value $\calJ_\ell(g_\infty,\bfi)\varphi(g_xg_\infty)$
at $g_\infty(\bfi)$.
Below, we construct cusp forms that generate $\pi_\infty(\ell)$ at the archimedean place and show that the corresponding modular forms do not vanish.

\section{Construction of cusp forms}

\subsection{Arthur's multiplicity formula}
We follow the notation and convention of \cite{Arthur_2013_book, 2015_Mok}.
For the discussion of Arthur parameters, take $A=\bbQ$ or $E$.
A discrete global Arthur parameter for $H_{A,n}$ is a formal sum
$\psi = \bigoplus_{i=1}^r \pi_i \boxtimes S_{d_i}$
where $S_{d_i}$ is the irreducible representation of $\SL_{2}(\bbC)$ of dimension $d_i$.
When $A=\Q$, each $\pi_i$ is a unitary selfdual cuspidal automorphic representation of $\GL_{N_i}(\bbA_\bbQ)$. 
The formal sum $\psi$ satisfy the following conditions.
\begin{itemize}
    \item $d_i$ is positive and $\sum_i N_i d_i = 2n+1$;
    \item If $d_i$ is odd, $\pi_i$ is orthogonal i.e., $L(s,\pi_i, \Sym^2)$ has a pole at $s=1$;
    \item If $d_i$ is even, $\pi_i$ is symplectic i.e., $L(s,\pi_i, \wedge^2)$ has a pole at $s=1$;
    \item The central character $\omega_i$ of $\pi_i$ satisfies $\omega_1^{d_1} \cdots \omega_{r}^{d_r} = \mathbf{1}$;
    \item If $i \neq j$ and $\pi_i \isom \pi_j$, then $d_i \neq d_j$.
\end{itemize}
When $A=E$, each $\pi_i$ is a unitary conjugate-selfdual cuspidal automorphic representation of $\GL_{N_i}(\bbA_E)$. 
The formal sum $\psi$ satisfy the following conditions.
\begin{itemize}
    \item $d_i$ is positive and $\sum_i N_i d_i = 2n$;
    \item If $d_i$ is odd, $\pi_i$ is conjugate-symplectic i.e., $L(s,\pi_i, \As^-)$ has a pole at $s=1$;
    \item If $d_i$ is even, $\pi_i$ is conjugate-orthogonal i.e., $L(s,\pi_i, \As^+)$ has a pole at $s=1$;
    \item If $i \neq j$ and $\pi_i \isom \pi_j$, then $d_i \neq d_j$.
\end{itemize}
Here $L(s,\pi_i,\As^+)$ and $L(s,\pi_i,\As^-)$ are the two
Asai $L$-functions for $E/\Q$, with the signs as in \cite{2015_Mok}.
For each global Arthur parameter $\psi$ and a place $v$ of $\bbQ$, we have a localization $\psi_v$.
The parameter $\psi_v$ gives rise to a finite multiset
$\Pi_{\psi_v}$ of irreducible representations of $H_{A,n}(\bbQ_v)$.
For the summand $\pi_i\boxtimes S_{d_i}$, let $e_{i,\bbA}$ be the corresponding symbol.
We have the (enhanced) component group
$\calS_\psi=\bigoplus_i(\Z/2\Z)e_{i,\bbA}$, with the central element $z_\psi=\sum_i e_{i,\bbA}$. At $v$ use the corresponding group
$\calS_{\psi_v}$ and central element $z_{\psi_v}$.
There is a canonical map $\Delta_v:\calS_\psi\to\calS_{\psi_v}$
with $\Delta_v(z_\psi)=z_{\psi_v}$.

Fix a global Whittaker datum. Its localization at $v$ determines
a map $\iota_v:\Pi_{\psi_v}\to\Irr(\calS_{\psi_v})$.
At each place choose $\pi_v\in\Pi_{\psi_v}$, unramified at almost all places,
and put $\eta_v=\iota_v(\pi_v)$. Let $\vep_\psi$ be the character of
$\calS_\psi/\langle z_\psi\rangle$ defined in \cite{Arthur_2013_book,2015_Mok}.
Suppose that $\eta_v(z_{\psi_v})=1$ for every $v$ and that $\eta_v=\mathbf1$
for almost all $v$. By Arthur's multiplicity formula, the restricted tensor product
$\bigotimes_v'\pi_v$ occurs in the discrete spectrum associated with $\psi$
if and only if $\prod_v(\eta_v\circ\Delta_v)=\vep_\psi$.

\subsection{Automorphic representations and their local components}
For a newform $f$, write $\pi_f$ for the unitary
cuspidal automorphic representation of $\GL_2(\bbA)$ associated with $f$.
For $A=\bbQ$ and even $n$, take $f\in S_{2k}(\SL_2(\bbZ))$ with
$k\equiv n/2\pmod2$ and put $\ell=k+n/2$. The global Arthur parameter of the
Ikeda lift is $\psi=\pi_f\boxtimes S_n\oplus\mathbf1\boxtimes S_1$.
For odd $n$, take $f\in S_{2k}(\SL_2(\bbZ))$ and $h\in S_{k+(n+1)/2}(\SL_2(\bbZ))$ with $k\equiv(n+1)/2\pmod2$. 
Put $\ell=k+(n+1)/2$. The global Arthur parameter of the Miyawaki lift is
$\psi=\pi_f\boxtimes S_{n-1}\oplus\Sym^2\pi_h\boxtimes S_1$, where
$\Sym^2\pi_h$ is the symmetric square lift of $\pi_h$
\cite[Theorem 9.3]{78_Gelbert_Jacquet}.

For $A=E$, let $\omega_{E/\Q}$ be the quadratic Hecke character of $E/\Q$, viewed also as a Dirichlet character with conductor $D$. 
We take a newform $f$ of weight $r>1$ without CM, that is, one that does not arise from a Hecke character of an imaginary quadratic field.
For odd $n$ take $f\in S_r(\SL_2(\Z))$, and for even $n$ take
$f\in S_r(\Gamma_0(D),\omega_{E/\Q})$. Put
$\psi =\pi_f^\BC\boxtimes S_n$.
Here $\BC$ denotes base change as in \cite{Langlands_Base_Change, 89_Arthur-Clozel}.
The base change $\pi_f^\BC$ is cuspidal and  conjugate self-dual with sign $(-1)^{n}$; see \cite[Lemma 11.1]{23_Ichino-Prasanna}. 
Put $\ell=n+r-1$.
The localization of $\psi$ at $\infty$ is given by
\[
 \psi_\infty=(\chi_{r-1}\boxtimes S_n)\oplus(\chi_{1-r}\boxtimes S_n),
 \qquad \chi_t(z)=(z/|z|)^t.
\]

For a newform $f\in S_{2k}(\SL_2(\bbZ))$, the component of $\pi_f$ at $p$
is $\pi_{f,p}=\mu_p\times\mu_p^{-1}$ for an unramified character $\mu_p$.
Let $a_p(f)$ denote the $p$-th Fourier coefficient of $f$ and let
$\alpha_p=\mu_p(p)$.
The set $\{\alpha_p,\alpha_p^{-1}\}$ are the roots of $X^2-a_p(f)p^{-k+1/2}X+1$ satisfying the Ramanujan bound $|\alpha_p|=1$.

\begin{lem}\label{lem_mu_p_nonquad}
If $f$ has integral Fourier coefficients, then $\mu_p^2\neq\mathbf1$ for every prime $p$.
\end{lem}
\begin{proof}
Fix a prime $p$. Since $\mu_p$ is unramified, it is determined by
$\alpha_p=\mu_p(p)$. If $\mu_p^2=\mathbf1$, then $\alpha_p=\alpha_p^{-1}=\pm1$.
Comparing the coefficient of $X$ in the Hecke polynomial above gives
\[
 a_p(f)=p^{k-1/2}(\alpha_p+\alpha_p^{-1})=\pm2p^{k-1}\sqrt{p}.
\]
This contradicts $a_p(f)\in\bbZ$.
\end{proof}

We use this lemma in weights $12$, $16$, and $18$ below.
We first consider the symplectic case.

\begin{lem}\label{lem:local-symplectic-packet}
Suppose $A=\bbQ$, $n$ is even and $f$ is as in Lemma~\ref{lem_mu_p_nonquad}. The packet $\Pi_{\psi_{p}}$ contains the irreducible unramified representation
$(\mu_p\circ\det_{\GL_n})\rtimes\mathbf1$.
\end{lem}
\begin{proof}
The localization $\psi_p$ is given by $(\mu_p\oplus\mu_p^{-1})\boxtimes S_n\oplus\mathbf1$. 
Since $\mu_p^2\neq\mathbf1$, the constituents $\mu_p \boxtimes S_n$ and $\mu_p^{-1} \boxtimes S_n$ are not self-dual.
By \cite[Proposition~8.11]{2017_Bin_Xu} and
\cite[Theorem~2.5]{202_Atobe_const_A-packet}, the only member of
$\Pi_{\psi_p}$ is $(\mu_p\circ\det_{\GL_n})\rtimes\mathbf1$, which is irreducible by \cite[Introduction, first theorem]{1992_Kudla-Rallis}.
\end{proof}

\begin{lem}\label{lem:local-miyawaki-atobe}
Let $A=\bbQ$, $n$ be odd, and write $\pi_{h,p}=\nu_p\times\nu_p^{-1}$. 
Suppose that $\mu_p^2\ne\mathbf1$.
Every irreducible constituent of
$I_{f,h,p}=(\mu_p\circ\det_{\GL_{n-1}})\times\nu_p^2\rtimes\mathbf1$
belongs to $\Pi_{\psi_{p}}$. 
For each constituent, the corresponding character $\eta_p$ satisfies $\eta_p \circ \Delta_p = \mathbf{1}$.

\end{lem}
\begin{proof}
Let $\sigma_p$ be a constituent of $\nu_p^2\rtimes\mathbf1$.
The character $\nu_p$ is unitary by the Ramanujan bound for $h$.
Thus $\sigma_p$ is tempered, and its $L$-parameter is
$\nu_p^2 \oplus \mathbf{1} \oplus \nu_p^{-2}$.
By \cite[Theorem 3.7]{2020_Atobe_Miyawaki} and the description of the local Arthur packet \cite[Proposition 8.11]{2017_Bin_Xu}, $(\mu_p\circ\det_{\GL_{n-1}})\rtimes\sigma_p$ is irreducible and belongs to $\Pi_{\psi_{p}}$.
The representation $(\mu_p\circ\det_{\GL_{n-1}})\rtimes\sigma_p$ is a constituent of $I_{f,h,p}$.
Conversely, every constituent of $I_{f,h,p}$ arises in this way.
Let $\eta_p$ be the corresponding character of $\calS_{\psi_{p}}$.
Since $\mu_p^2\neq\mathbf{1}$, the symbol corresponding to $\pi_f\boxtimes S_{n-1}$ has zero image under $\Delta_p$.
The other generator corresponding to $\Sym^2 \pi_h \boxtimes S_1$ therefore maps to $z_{\psi_{p}}$.
Since $\eta_p(z_{\psi_{p}}) = 1$, we have $\eta_p \circ \Delta_p = \mathbf{1}$.
\end{proof}

For the case $A=E$, at a nonsplit prime $p$ the newform $f$ used below satisfies
\[
\pi_{f,p}=\begin{cases}
\mu_p\times\mu_p^{-1} & n\text{ odd},\\
\mu_p\times\mu_p^{-1}\omega_{E_p/\bbQ_p} & n\text{ even},
\end{cases}
\qquad \mu_p^\BC=\mu_p\circ N_{E_p/\bbQ_p},
\]
where $\mu_p$ is unramified and unitary; see \cite[Appendix A]{2025_HMY_Kodaira_dimension_ball_quotients}. Thus the $L$-parameter of $\pi_{f,p}^\BC$ is $\mu_p^\BC\oplus(\mu_p^\BC)^{-1}$.
By \cite[Theorem 1.1]{1997_Kudla-Sweet}, at nonsplit primes, the representation 
$I_{f,p} \defeq (\mu_p^\BC\circ\det_{\GL_n})\rtimes\mathbf1$
is reducible if and only if
\[
\mu_p^2 = \mu_p^\BC|_{\mathbb Q_p^\times}
=
\begin{cases}
\omega_{E_p/\bbQ_p} & n\text{ odd},\\
\mathbf1 & n\text{ even}.
\end{cases}
\]

\begin{lem}\label{lem_unitary_local_noninert}
Let $A=E$. If $p$ splits in $E$, the packet $\Pi_{\psi_p}$ contains an
unramified representation of $\GL_{2n}(\Q_p)$. If $p$ is nonsplit, every
irreducible constituent of $I_{f,p}$ belongs to $\Pi_{\psi_p}$.
\end{lem}
\begin{proof}
At a split prime, the parameter for $\GL_{2n}(\Q_p)$ gives an unramified
representation. At a nonsplit prime where $I_{f,p}$ is irreducible, the
assertion follows from \cite[Proposition~8.4.1]{2015_Mok}.

Suppose now that $I_{f,p}$ is reducible. The reducibility criterion above
gives $(\mu_p^\BC)^2=\mathbf1$, since both $\mathbf1$ and
$\omega_{E_p/\Q_p}$ are trivial on norms from $E_p^\times$.
Interchanging the two $\SL_2$ factors gives the tempered parameter
\[
 \widehat{\psi_p}=(\mu_p^\BC\boxtimes S_n\boxtimes S_1)^{\oplus 2}.
\]
Its packet consists of the irreducible constituents of
$\mu_p^\BC\mathrm{St}_{\GL_n}\rtimes\mathbf1$, where
$\mathrm{St}_{\GL_n}$ is the Steinberg representation of $\GL_n(E_p)$.
                By \cite[Theorem~5.4.1]{AGIKMS_LIR_2024}, the Aubert duals of these
constituents form $\Pi_{\psi_p}$. Aubert duality commutes with normalized
parabolic induction up to sign \cite[Theorem~1.7]{1995_Aubert} and sends
$\mu_p^\BC\mathrm{St}_{\GL_n}$ to $\mu_p^\BC\circ\det_{\GL_n}$.
Thus every irreducible constituent of $I_{f,p}$ belongs to $\Pi_{\psi_p}$.
\end{proof}

\subsection{Fixed vectors and automorphy}\label{subsec:fixed-automorphy}

Take $K_\Gamma$ as in Section~\ref{subsec:classical-symplectic-unitary}.
For $A=\bbQ$ in genus $15$, take $k=8$.
In the other genera take $k = 6$ if $n\equiv0,3\pmod4$
and $k = 9$ if $n\equiv1,2\pmod4$. 
Each of the corresponding newforms $f$ has integral Fourier coefficients.
For odd $n$ choose a newform $h$ of weight
$\ell=k+(n+1)/2$ for $\SL_2(\Z)$. 
Such a cusp form exists in the ranges of Table~\ref{tab:cusp-weights} for $A=\Q$. 
For $n=11,15,17$ we have $\ell=2k=n+1$ and take $h=f$.

For $H_{E,n}$ take $f=\Delta$ if $n$ is odd. For even $n$,
Theorem~\ref{thm_existence_cusp_rep_no_twist} gives a newform $f$ without CM of weight
$r=c_E+1$, level $D$, and character $\omega_{E/\Q}$.

For $A=\Q,E$ take $\pi_\infty=\pi_\infty(\ell)$ and put
$\eta_\infty=\iota_\infty(\pi_\infty)$. 
For $A=\Q$, or for $A=E$ and $p$ nonsplit, the induced representations $\mu_p \circ \det_{\GL_n} \rtimes \mathbf{1}, I_{f,p}$ and $I_{f,h,p}$ fulfill the condition as in Lemma \ref{lem:fixed}.
Hence, the representations contains a $K_{\Gamma,p}$-fixed vector.
Let $\pi_p$ be an irreducible representation in the packet containing such a vector.
When $p$ splits in $E$, take the unramified
member $\pi_p$. The group $K_{\Gamma,p}$ preserves a lattice in $\Q_p^{2n}$
and thus lies in a conjugate of $\GL_{2n}(\Z_p)$, so $\pi_p^{K_{\Gamma,p}}\ne0$.
Choose such a vector at each finite prime, taking the spherical vector in the unramified representation at almost all primes, and put
$\eta_p=\iota_p(\pi_p)$.

For $A=B$, we consider Yamana's Ikeda type lifting \cite{Yamana} for quaternionic unitary groups.
In this case, the lifting of $\Delta \in S_{12}(\SL_2(\bbZ))$ is of weight $\ell = n+5$.
Let $\pi$ be the cuspidal automorphic representation associated with the lifting.

\begin{lem}\label{lem_automorphy}
The representation $\pi=\bigotimes_{v}'\pi_v$ is automorphic and discrete.
\end{lem}
\begin{proof}
The cases $A=B$ and $A=\bbQ$ with $n$ even follow from the definition.
In the unitary case,  $\calS_\psi=\langle z_\psi\rangle$ and $\Pi_{\psi_\infty}$ contains $\pi_\infty(\ell)$ by \cite[Corollary 4.8]{Horinaga_LW_unitary}.
Since $\Delta_v(z_\psi)=z_{\psi_v}$ and $\eta_v(z_{\psi_v})= \mathbf{1}$ at every place,
we have $\prod_v(\eta_v\circ\Delta_v)=\mathbf{1}=\vep_\psi$ for every choice of the
representations $\pi_v$ at the finite places.
By Arthur's multiplicity formula, $\pi$ is automorphic and discrete.

It remains to treat $A=\bbQ$ with odd $n$. 
Use the Whittaker datum $\mathrm{Wh}_+$ of \cite{MR_ht_wt}. 
Suppose that $\ell \leq n$.
With $\rho_t=\Ind_{W_\C}^{W_\R}(z/|z|)^t$, where $W_\C$ and $W_\R$ are the Weil groups of $\C$ and $\R$, the parameter at infinity is
\[
 \psi_\infty=\rho_{2k+n-1}\boxtimes S_1\oplus\rho_{2k-1}\boxtimes S_{n-1}\oplus\sgn\boxtimes S_1.
\]
By \cite[Theorem 7.1]{MR_ht_wt}, the representation $\pi_\infty$ is contained in the packet $\Pi_{\psi_\infty}$.
For the lowest weight $\ell=k+(n+1)/2 \leq n$, \cite[Proposition~18.3]{MR_ht_wt} gives
\[
\eta_\infty(e_{\rho_{2k+n-1}\boxtimes S_1})=1,\quad
\eta_\infty(e_{\rho_{2k-1}\boxtimes S_{n-1}})
=\eta_\infty(e_{\sgn\boxtimes S_1})=(-1)^{(n-1)/2}.
\]
Since $\pi_f$ and $\pi_h$ are unramified at every finite place, the root number is given by
\[
\vep_\psi(e_{\Sym^2\pi_h\boxtimes S_1})
=\vep_\psi(e_{\pi_f\boxtimes S_{n-1}})=(-1)^{k+1}
\]
by \cite[Section~8.2]{2019_Chenevier_Lannes}.
Under localization these generators map to
$e_{\rho_{2k+n-1}\boxtimes S_1}+e_{\sgn\boxtimes S_1}$ and
$e_{\rho_{2k-1}\boxtimes S_{n-1}}$, respectively. Since $k\equiv(n+1)/2\pmod2$, the character $\eta_\infty \circ \Delta_\infty$ equals $\vep_\psi$. 
The characters $\eta_v$ for $v<\infty$ satisfy $\eta_v \circ \Delta_v = \mathbf{1}$ by Lemma~\ref{lem:local-miyawaki-atobe}. Arthur's multiplicity formula proves the assertion.
The case $\ell = n+1$ follows from the same method and \cite[Lemma 9.4 and Appendix A]{2015_Chenevier_Renard}.\end{proof}

\subsection{Cuspidality and connected components}
\begin{lem}\label{lem_cuspidality}
The representation $\pi$ constructed above is cuspidal.
\end{lem}
\begin{proof}
The case $A=B$ follows from the definition.
Below we assume $A=\bbQ$ or $E$.
If $\pi$ were in the residual spectrum, by \cite[Section~1.3]{2008_Moeg_noncusp}, $\pi_\infty$ is a
subquotient of $\sigma_\infty|\det|^{(d-1)/2}\rtimes\pi'_\infty$,
where $(\sigma,d)$ is $(\pi_f,n)$ for $A=\bbQ$ with even $n$,
$(\pi_f,n-1)$ for $A = \bbQ$ with odd $n$, and $(\pi_f^\BC,n)$ for $A=E$.
Here $\pi'_\infty$ belongs to the packet for the parameter obtained by
replacing $S_d$ with $S_{d-2}$ and omitting the summand when $d=2$.

For $A=\bbQ$, the maximal compact subgroup of this Levi subgroup is $\O(2)\times\U(n-2)$. 
The scalar $K_\infty$-type $\tau_\ell$ restricts trivially to $\SO(2)$, whereas the discrete series $\pi_{f,\infty}$ of weight $2k$ has no $\SO(2)$-fixed vector. 
By Frobenius reciprocity, we have
\[
\Hom_{\U(n)}(\tau_\ell, \sigma_\infty|\det|^{(d-1)/2} \rtimes \pi_\infty') \isom 
\Hom_{\SO(2) \times \U(n-2)}(\tau_\ell|_{\SO(2) \times \U(n-2)}, \sigma_\infty|_{\SO(2)} \boxtimes \pi_\infty'|_{\U(n-2)}) = 0.
\]
Hence $\pi_\infty(\ell)$ cannot occur in a subquotient.

For $A=E$, the Levi subgroup is $\GL_2(\mathbb C)\times\U(n-2,n-2)$. 
The maximal compact subgroup $\U(2)$ of $\GL_2(\bbC)$ embeds diagonally in the two factors of $K_\infty$, and $\tau_\ell$ restricts trivially to it. 
Since $\pi_{f,\infty}^\BC=\chi_{r-1}\times\chi_{1-r}$ and $r>1$, the assertion follows from the same method as the case $A=\bbQ$.
\end{proof}

\begin{lem}\label{lem:unitary-no-self-twist}
Let $\tau$ be a cuspidal automorphic representation of $\GL_2(\bbA_\bbQ)$.
If $\tau$ is unramified at every finite place, there is no nontrivial Hecke character $\xi$ such that $\tau \isom \tau \otimes (\xi \circ \det_{\GL_2})$.
If no nontrivial Hecke character $\xi$ of $\bbA_\bbQ^\times$ satisfies $\tau \isom \tau \otimes (\xi \circ \det_{\GL_2})$, then no nontrivial Hecke character $\eta$ of $\bbA_E^\times$ satisfies $\tau^\BC \isom \tau^\BC \otimes (\eta \circ \det_{\GL_2})$. 

\end{lem}
\begin{proof}
The first case follows from \cite[Lemma~11.7]{Langlands_Base_Change}.
For the second, by the assumption and \cite[Theorem~9.3]{78_Gelbert_Jacquet}, $\Sym^2\tau$ is cuspidal.
If its base change to $E$ were not cuspidal, then $\Sym^2\tau \otimes (\omega_{E/\bbQ} \circ \det_{\GL_3}) \isom \Sym^2\tau$ by \cite[Chapter~3, Theorem~4.2]{89_Arthur-Clozel}.
Its central character would then give $\omega_{E/\Q}^3=1$, which is impossible.
As the symmetric square lifting commutes with the base change lifting, $\Sym^2(\tau^\BC)$ is
cuspidal. 
Hence, $\tau^\BC$ has no nontrivial self twist by \cite[Remark~9.9]{78_Gelbert_Jacquet}.
\end{proof}

\begin{proof}[Proof of Theorem~\ref{thm:intro-cusp-forms}]
By Lemmas~\ref{lem_automorphy} and~\ref{lem_cuspidality}, $\pi$ is cuspidal automorphic.
For $A=\Q,E$, tensor the fixed vectors chosen in
Section~\ref{subsec:fixed-automorphy} with a vector of lowest weight at infinity.
This gives a nonzero cusp form $\varphi$ fixed by $K_\Gamma$. Its vector at infinity
has $K_\infty$-type $\tau_\ell$ and is annihilated by $\frakp_-$, so $\varphi$ is
holomorphic of weight $\ell$. For $A=\Q$ there is only one connected component by strong approximation.
Hence, the assertion for $A=\bbQ$ holds.

For $A=E$, let $\xi$ be a nontrivial character of $\Cl_{K_\Gamma}$.
By Hilbert's Theorem~90, $\eta_\xi(x)=\xi(x/\bar x)$ is a nontrivial
Hecke character of $E$.
Suppose that $\pi\simeq\pi\otimes(\xi\circ\det)$. 
Let $w$ be a finite place of $E$ above $p$ such that all data are unramified and $\{\alpha_w, \alpha_w^{-1}\}$ be the set of Satake parameters of $\pi_{f,w}^\BC$.
Let $\varpi_w$ be a uniformizer in $\calO_{F,w}$.
The set of Satake parameters for $\pi$ (resp.~$\pi\otimes(\xi\circ\det)$) at $p$ is given by $\{\alpha_w p^{(n+1)/2-j}, \alpha_w^{-1} p^{(n+1)/2-j} \mid j = 1,\ldots,n\}$ (resp.~$\{\alpha_w p^{(n+1)/2-j}\eta_{\xi,w}(\varpi_w), \alpha_w^{-1} p^{(n+1)/2-j}\eta_{\xi,w}(\varpi_w) \mid j = 1,\ldots,n\}$).
Since $|\alpha_w|=|\eta_{\xi,w}(\varpi_w)| = 1$, we have $\pi_{f,w}^\BC \isom \pi_{f,w}^\BC \otimes \eta_{\xi,w}$ by $\{\alpha_w,\alpha_w^{-1}\} = \{\alpha_w\eta_{\xi,w}(\varpi_w),\alpha_w^{-1}\eta_{\xi,w}(\varpi_w)\}$.
By the strong multiplicity one \cite[Theorem 4.4]{1981_Jacquet_Shalika_gen_mult_one}, $\pi_{f}^\BC$ is isomorphic to $\pi_f^\BC \otimes(\eta_\xi\circ\det)$.
Using Lemma \ref{lem:unitary-no-self-twist} and our assumption on $f$, we conclude that $\eta_\xi = \mathbf{1}$, which is a contradiction.
Hence, Lemma \ref{prop:common-components} shows that the restrictions of the cusp form to each connected component are nonvanishing.

For $A=B$, by \cite[p.2495]{Yamana}, the representations $\pi_p$ for finite primes $p$ are isomorphic to $|\nu|_p^{s_p} \rtimes \mathbf{1}$ for some $\mathrm{Re}(s_p)=0$.
Since $\pi_p$ are irreducible as in \cite[Proposition 4.1 (1)]{Yamana}, there exists a nonzero $K_{\Gamma,p}$-fixed vector $v_p$ in $\pi_p$.
Taking a lowest weight vector $v_\infty$ at the infinite place, $\otimes_w v_w$ gives rise to a $K_\Gamma$-fixed cusp form $\varphi$ on $H_{B,n}(\bbA)$ that is of weight $\ell = n+5$ and annihilated by $\frakp_-$.
It remains to show that every restriction of $\varphi$ to each connected component is nonvanishing.
Let $\chi$ be a character of $\Cl_{K_\Gamma}$ and $\widetilde{\chi}$ be the character of $H_{B,n}(\bbA)$ defined by $\chi \circ \sn$.
Suppose that $\pi \isom \pi \otimes \widetilde{\chi}$.
Let $\{\alpha_p, \alpha_p^{-1}\}$ be the Satake parameters for $\pi_{\Delta,p}$.
Then, the set of Satake parameters of $\pi$ (resp.~$\pi \otimes \widetilde{\chi}$) is given by $\{\alpha_pp^{n+1/2-j}, \alpha_p^{-1}p^{n+1/2-j} \mid j = 1,\ldots,2n\}$ (resp.~$\{\chi_p(p)\alpha_pp^{n+1/2-j}, \chi_p(p)\alpha_p^{-1}p^{n+1/2-j} \mid j = 1,\ldots,2n\}$).
Comparing the Satake parameters, by $|\alpha_p| = |\alpha_p^{-1}| = 1$, we have $\{\alpha_p,\alpha_p^{-1}\} = \{\alpha_p\chi_p(p), \alpha_p^{-1}\chi_p(p)\}$, so $\pi_{\Delta,p} \isom \pi_{\Delta, p} \otimes (\chi_p \circ \det)$ for almost all $p$.
By the strong multiplicity one, we have $\pi_{\Delta} \isom \pi_{\Delta} \otimes (\chi \circ \det)$.
By Lemma \ref{lem:unitary-no-self-twist}, such $\chi$ is trivial.
Thus $\varphi$ is nonzero on each connected component by Lemma~\ref{prop:common-components}.
\end{proof}

\section{Singularities and Kodaira dimension}
\label{sec:boundary}

\subsection{Singularities of arithmetic quotients}
For the description of the toroidal charts we follow Ash, Mumford, Rapoport, and Tai
\cite[Chapter~III, Section~4.3 and p.~158]{ash2010smooth}.
Let $\DD$ be a Hermitian symmetric domain and let $G$
be a connected semisimple adjoint group over $\Q$ with
$G(\R)^\circ=\Aut(\DD)^\circ$. Let
$\GG\subset G(\Q)\cap G(\R)^\circ$ be arithmetic. 

For a rational boundary component $F$ of $\calD$, let $P_F$ be its stabilizer in $G$ and $U_F$ the center of the unipotent radical of $P_F$. Put $\GG_F=\GG\cap P_F(\Q)$ and $\Lambda_F=\GG\cap U_F(\Q)$.
The group $U_F$ is a vector group, and $\Lambda_F$ is a full lattice in
$U_F(\R)$, that is, it is discrete and the natural map
$\Lambda_F\otimes_\Z\R\to U_F(\R)$ is an isomorphism.
Let $C_F\subset U_F(\R)$ be the open cone in the realization of $\DD$ as a Siegel domain.
Put $\DD(F)=U_F(\C)\cdot\DD$ in the
compact dual and $B_F=\DD(F)/U_F(\C)$. The projection
$\pi_F:\DD(F)\to B_F$ is a principal $U_F(\C)$-bundle. For each
$b\in B_F$, the set $\DD_b=\DD\cap\pi_F^{-1}(b)$ is a nonempty convex
open subset of the affine fiber.

We choose a projective $\GG$-admissible collection of fans $\Sigma=(\Sigma_F)$
such that each $\Sigma_F$ is regular with respect to $\Lambda_F$ and the
stabilizer in $\GG_F$ of every cone fixes it pointwise.
Start with a projective admissible collection for a neat normal subgroup of
finite index in $\GG$ \cite[Chapter~IV, Section~2]{ash2010smooth}. The sum of
the support functions of its finitely many $\GG$-translates defines their
common refinement. Its barycentric subdivision is again projective, since
it is obtained by star subdivisions. Each cone of this subdivision corresponds
to a flag of faces of distinct dimensions. Its stabilizer therefore preserves
each face in the flag and fixes the primitive generator of each ray of the cone.
Now refine one representative of each orbit of cones in order of increasing
dimension, extending the subdivisions on the faces, by
\cite[Lemma~5.22 and Proposition~5.23]{pink1989arithmetic}.
Since each cone is fixed pointwise by its stabilizer, these choices extend
$\GG$-equivariantly and give a projective refinement regular for $\Lambda_F$.
The stabilizer of each new cone preserves the smallest cone of the barycentric
subdivision containing it and hence fixes the new cone pointwise.

Let $\mathcal X_F=(\DD/\Lambda_F)^{\Sigma_F}$ be the partial toroidal
compactification in the torus embedding of $\DD(F)/\Lambda_F$ over $B_F$.
The primitive ray generators of each cone extend to a basis of $\Lambda_F$.
Thus its toric chart is smooth, and $\mathcal X_F$ is smooth over $B_F$.
The stratum over $F$ corresponds to the cones whose relative interiors lie in $C_F$. 
Near the image of a point $x$ over this stratum, the global compactification is the quotient of a neighborhood of $x$ in $\mathcal X_F$ by its finite stabilizer $H_x\subset\GG_F/\Lambda_F$. 
See \cite[Proposition~III.6.10]{ash2010smooth}.

If a linear transformation $g$ of finite order on a complex vector space $V$ has eigenvalues $\exp(2\pi\sqrt{-1}\alpha_j)$ with $0\leq\alpha_j<1$, we define its \emph{age} by $\age(g;V)=\sum_j\alpha_j$.
Recall that $X$ has canonical singularities if it
is normal, $K_X$ is $\Q$-Cartier, and $K_{\widetilde X}-f^*K_X$ is effective
on a resolution $f:\widetilde X\to X$. 
By a quasireflection, we mean an action that fixes a hyperplane pointwise and acts nontrivially on its normal line. 
For a finite group of linear transformations without quasireflections, the Reid--Shepherd-Barron--Tai criterion characterizes canonical singularities by the inequality
$\age(g)\geq1$ for every nonidentity element \cite[Section~4]{reid1987young}.
We linearize the neighborhood of a point and apply the criterion to the action of the stabilizer on the tangent space.

\begin{prop}\label{prop:boundary-age}
Let $x\in\mathcal X_F$ lie over the stratum associated with $F$ and
let $1\ne h\in H_x$. There is a nonidentity element of finite order
$s\in G(\Q)\cap G(\R)^\circ$ and a point $p\in\DD$ fixed by $s$
such that
\[
 \age(h;T_x\mathcal X_F)\geq\age(s;T_p\DD).
\]
\end{prop}

\begin{proof}
Let $r$ be the order of $h$ and choose a representative $\gamma\in\GG_F$.
Since $\gamma^r\in\Lambda_F$ and $U_F$ is a vector group, $\gamma^r$ has a unique $r$th root $u$ in $U_F(\Q)$. Conjugation by $\gamma$ fixes $u$ by uniqueness. 
Thus $s=\gamma u^{-1}$ belongs to $P_F(\Q)\cap G(\R)^\circ$ and satisfies $s^r=1$. 
If $s=1$, then $\gamma\in\GG\cap U_F(\Q)=\Lambda_F$, contrary to $h\ne1$.

Let $b\in B_F$ be the image of $x$. Since $h$ fixes $x$, its representative
$\gamma$ fixes $b$. The translation $u$ acts trivially on $B_F$, so
$s=\gamma u^{-1}$ fixes $b$ and induces the same action as $h$ on $T_bB_F$.
The element $s$ acts affinely on $\pi_F^{-1}(b)$ and preserves its convex
open subset $\DD_b$. For any $p_0\in\DD_b$, the barycenter
$p=r^{-1}\sum_{j=0}^{r-1}s^jp_0$ lies in $\DD_b$ and satisfies $sp=p$.

Write $\mathcal X_{F,b}$ for the fiber of $\mathcal X_F\to B_F$ over $b$.
The projections to $B_F$ give the exact sequences
\[
0\longrightarrow U_F(\C)\longrightarrow T_p\DD
  \longrightarrow T_bB_F\longrightarrow 0,\qquad 
0\longrightarrow T_x\mathcal X_{F,b}\longrightarrow T_x\mathcal X_F
  \longrightarrow T_bB_F\longrightarrow 0.
\]
The first sequence is equivariant for $s$, the second for $h$, and the actions
on $T_bB_F$ agree. On $U_F(\C)$ the element $s$ acts by conjugation.
Write $\phi$ for this action, which coincides with conjugation by $\gamma$
since $U_F$ is abelian.
Let $O_\sigma\subset\mathcal X_{F,b}$ be the torus orbit containing $x$
and put $N_\sigma=\Lambda_F\cap\operatorname{span}_{\R}(\sigma)$.
On $O_\sigma$ the element $h$ acts by a translation composed with the
automorphism induced by $\phi$
\cite[proof of Lemma~6.3]{2026_Horinaga-Maeda_sing_unitary}.
Since $h$ fixes $x$, its derivative on
$T_xO_\sigma\simeq(\Lambda_F/N_\sigma)\otimes_\Z\C$ is induced by $\phi$.
The choice of fans ensures that $\phi$ fixes $\sigma$ pointwise and hence
acts trivially on $N_\sigma\otimes_\Z\C$. Thus
$\age(h;T_xO_\sigma)=\age(s;U_F(\C))$.
By additivity of age in the two sequences above and in the tangent sequence
of $O_\sigma\subset\mathcal X_{F,b}$, we obtain
\[
\age(h;T_x\mathcal X_F)-\age(s;T_p\DD)
   =\age(h;T_x\mathcal X_{F,b}/T_xO_\sigma)\geq 0. \qedhere
\]
\end{proof}

\begin{thm}\label{thm:boundary-transfer}
Suppose that every nonidentity element of finite order
$s\in G(\Q)\cap G(\R)^\circ$ fixing a point $p\in\DD$ satisfies
$\age(s;T_p\DD)\geq1$.
Then $\GG\backslash\DD$ admits a projective toroidal compactification with canonical singularities and has no branch divisors.
\end{thm}

\begin{proof}
Let $x \in \calX_F$.
Every nonidentity element of $H_x$ has age at least one, by hypothesis in the interior and by Proposition~\ref{prop:boundary-age} at the boundary. 
Thus there are no quasireflections, so there are no branch divisors and the Reid--Shepherd-Barron--Tai criterion gives canonical singularities at $x$.
For the absence of ramification along the boundary divisors, see also \cite[Theorem~1]{ma2022boundary}.
\end{proof}

\begin{prop}\label{prop:central-descent}
Let $G$, $\DD$, and an arithmetic subgroup
$\Gamma\subset G(\Q)\cap G(\R)^\circ$ satisfy the hypotheses of
Theorem~\ref{thm:boundary-transfer}. Let $q:\widetilde G\to G$ be a central
isogeny of connected algebraic groups over $\Q$. Put
$\Gamma_0=\Gamma\cap q(\widetilde G(\Q))$.
Then $\Gamma_0$ is normal of finite index in $\Gamma$, and
$\kappa(\Gamma_0\backslash\DD)=\kappa(\Gamma\backslash\DD)$.
\end{prop}

\begin{proof}
Set $Z=\ker q$. The connecting homomorphism $\Gamma\to H^1(\Q,Z)$
has kernel $\Gamma_0$. The arithmetic group $\Gamma$ is finitely generated
by \cite[Theorem~6.12]{BorelHarishChandra1962}, since finite generation is
invariant under commensurability. Its image in $H^1(\Q,Z)$ is therefore
a finitely generated abelian group killed by the exponent of $Z$, and is finite.

The isogeny $q$ identifies the corresponding unipotent radicals over
$\Q$, so $\Gamma_0\cap U_F(\Q)=\Gamma\cap U_F(\Q)$ at every cusp.
With the same fans for both groups, the induced finite morphism between
the toroidal compactifications is locally a map between quotients of the
same smooth chart. By Theorem~\ref{thm:boundary-transfer}, both
compactifications have canonical singularities and the chart maps have
no branch divisors. The induced morphism is therefore unramified in codimension
one, so the canonical divisor of the source is the pullback of the canonical
divisor of the target. Since finite surjective pullback preserves Iitaka
dimension, the two quotients have the same Kodaira dimension.
\end{proof}

\subsection{Singularities in type DIII and the proof of Theorem~\ref{mainthm:kodaira}}
Take $\DD=\DD^{\mathrm{DIII}}_m$ in Section~\ref{sec:boundary} and assume
that $G$ is absolutely simple. 
At a point $p$ of $\DD^{\mathrm{DIII}}_m$, the stabilizer of $p$ is isomorphic to $\mathrm U(m)/\{\pm I\}$. 
Choose the standard representation of $\mathrm U(m)$ so that its exterior square is isomorphic to the holomorphic tangent space as a representation of $\U(m)$. 
If $e_1,\ldots,e_m$ are the characters of the diagonal torus, the weights on the tangent space are $e_i+e_j$ for $i<j$, and the roots of the adjoint representation are $\pm e_i\pm e_j$.

\begin{lem}\label{lem:cyclotomic}
Let $T\in\GL(V_{\Q})$ have finite order. If every eigenvalue of $T$
has real part greater than $1/2$, then $T=I$.
\end{lem}

\begin{proof}
If an eigenvalue has order $d>1$, every primitive $d$-th root of unity
occurs among the eigenvalues, since the characteristic polynomial has
rational coefficients. Put $\zeta_d=\exp(2\pi\sqrt{-1}/d)$. 
We have $\varphi(d)^{-1}\sum_{a\in(\Z/d\Z)^\times}\zeta_d^a
=\mu(d)/\varphi(d)\leq1/2$, where $\mu$ is the M\"obius function and
$\varphi$ is Euler's totient function.
This contradicts the hypothesis on their real parts. 
\end{proof}

\begin{prop}\label{prop:interior-age}
Let $s\in G(\Q)\cap G(\R)^\circ$ be a nonidentity element of finite
order. If $s$ fixes $p\in\DD^{\mathrm{DIII}}_m$ and $m\geq7$, then
$\age(s;T_p\DD^{\mathrm{DIII}}_m)\geq1$.
\end{prop}

\begin{proof}
Lift $s$ to $\mathrm U(m)$ and diagonalize the lift, with eigenvalues
$\lambda_1,\ldots,\lambda_m$. Its eigenvalues on $T_p\DD^{\mathrm{DIII}}_m$
are $\lambda_i\lambda_j$ for $i<j$. Write
$\lambda_i\lambda_j=\exp(2\pi\sqrt{-1}\alpha_{ij})$, where
$0\leq\alpha_{ij}<1$, and set $\alpha_{ji}=\alpha_{ij}$.
Suppose that $\age(s)=\sum_{i<j}\alpha_{ij}<1$. For four distinct indices,
the sums $\alpha_{ik}+\alpha_{jl}$ and $\alpha_{ij}+\alpha_{kl}$ are
congruent modulo $\Z$. Both lie in $[0,\age(s)] \subset [0,1)$, so they are equal.
Fix $i\ne j$ and put $C=\{1,\ldots,m\}\setminus\{i,j\}$.
Sum $\alpha_{ik}+\alpha_{jl}\geq\alpha_{ij}$ over the ordered pairs of
distinct $k,l\in C$. 
Hence, $(m-3)\sum_{k \in C} (\alpha_{i,k} + \alpha_{j,k})=\sum_{k,\ell \in C, k \neq \ell} (\alpha_{i,k} + \alpha_{j,\ell}) \geq (m-2)(m-3)\alpha_{i,j}$ and $\sum_{k \in C} (\alpha_{i,k} + \alpha_{j,k}) \geq (m-2)\alpha_{i,j}$.
Since $\age(s)\geq\alpha_{ij}+\sum_{k\in C}(\alpha_{ik}+\alpha_{jk})$, we obtain
\begin{equation}\label{eq:angle-small}
 0\leq\alpha_{ij}\leq\frac{\age(s)}{m-1}<\frac{1}{m-1}\leq\frac16.
\end{equation}
Consider the adjoint representation. Its zero weights give eigenvalue one, and its
roots $\pm e_i\pm e_j$ give the products
$\lambda_i^{\pm 1}\lambda_j^{\pm 1}$.
For $k$ distinct from $i,j$, we have
\[
 \frac{\lambda_i}{\lambda_j}
 =\frac{\lambda_i\lambda_k}{\lambda_j\lambda_k}
 =\exp\bigl(2\pi\sqrt{-1}(\alpha_{ik}-\alpha_{jk})\bigr).
\]
By \eqref{eq:angle-small}, every eigenvalue $\lambda_i^{\pm 1}\lambda_j^{\pm 1}$ has an
argument of absolute value strictly less than $2\pi/6$. 
Its real part is therefore greater than $1/2$.
Since $s\in G(\Q)$ has finite order, Lemma~\ref{lem:cyclotomic} gives
$\Ad(s)=I$. The adjoint representation is faithful, so $s=1$, a contradiction.
\end{proof}

\begin{thm}\label{thm:DIII-compactification}
For $m\geq7$, the quotient $\GG\backslash\DD^{\mathrm{DIII}}_m$ admits a
projective toroidal compactification with canonical singularities and has no branch divisors.
\end{thm}
\begin{proof}
The assertion follows from Theorem \ref{thm:boundary-transfer} and Proposition \ref{prop:interior-age}.
\end{proof}
The following example shows the sharpness of the above bound.
\begin{ex}\label{ex:sharpness}
In the definite quaternion algebra $B'=(-1,-3)_\Q$, take the standard
generator $j$ with $j^2=-3$. The element $u=(1+j)/2$ has reduced norm
one and order six, since $u^2-u+1=0$. On $(B')^m$ consider the
skew Hermitian form $(v,w)\mapsto v^*(jI_m)w$. 
Its group of real isometries is $\SO^*(2m)$. With $\R[j]\simeq\C$, its
maximal compact subgroup at the base point is
$K=\{A\in\GL_m(\R[j])\mid\bar A^{\,t}A=I_m\}\simeq\mathrm U(m)$.
Indeed, the entries of $A\in K$ commute with $j$, so $A^*(jI_m)A=jA^*A=jI_m$.
Since $u\bar u=1$, the rational isometry $g=\diag(u,1,\ldots,1)$ lies
in $K$ and fixes the base point.
After replacing $g$ by its inverse if necessary, we find that its eigenvalues
on the tangent space are
$\exp(2\pi\sqrt{-1}/6)$ with multiplicity $m-1$ and $1$ with
multiplicity $(m-1)(m-2)/2$.
In particular, its age is $(m-1)/6$.

Let $s$ be the image of $g$ in the adjoint group $G$. It still has order
six. 
For any full lattice $L$ in $\operatorname{Lie}(G)(\Q)$, the lattice
$\sum_{a=0}^{5}\Ad(s)^aL$ is preserved by $s$. Its stabilizer
$\Gamma_1$ in $G(\Q)\cap G(\R)^\circ$ is arithmetic and contains $s$.
As in \cite[proof of Theorem~9.23]{2026_Horinaga-Maeda_sing_unitary}, take a torsion free normal subgroup $\Gamma_0$ of finite index in $\Gamma_1$
and put $\Gamma=\langle s\rangle\Gamma_0$. The stabilizer of the base
point in $\Gamma_0$ is trivial, so its stabilizer in $\Gamma$ is $\langle s\rangle$.
The fixed space of each nonidentity power of $s$ on the tangent space
has codimension $m-1$, so there are no quasireflections for $m\geq3$.
For $3\leq m\leq6$ the age $(m-1)/6$ is less than one, so the quotient
singularity is not canonical.
\end{ex}

We return to the arithmetic subgroups $\Gamma\subset H_{A,n}(\Q)$ of
Section~\ref{subsec:intro-kodaira} for general $A$. 
Let $C_0 = \{\pm I_{2n}\}$ for $A=\bbQ,B$ and $C_0 = \{z \cdot I_{2n} \mid z \in \calO^\times_E\}$ for $A=E$.
These scalars act trivially on the domain. Thus $\Gamma$ and $\Gamma C_0$
have the same image in $\Aut(\frakH_{A,n})$ and, for the same fans,
the same coarse toroidal compactification. In particular, $X_{\Gamma C_0}=X_\Gamma$.
Combining \cite{tai1982kodaira, 2026_Horinaga-Maeda_sing_unitary} with the preceding theorem for arithmetic subgroups $\Gamma$ on $\DD^{\mathrm{DIII}}_{2n}$ gives the following corollary.

\begin{cor}[{\cite[Propositions~4.2 and~5.3]{tai1982kodaira},
\cite[Theorem~2.2]{2026_Horinaga-Maeda_sing_unitary}, and
Theorem~\ref{thm:DIII-compactification}}]\label{cor:canonical-toroidal}
Let $\Gamma\subset H_{A,n}(\Q)$ be arithmetic.
The quotient $X_\Gamma$ admits a projective toroidal compactification $\overline X_\Gamma^\Sigma$ with canonical singularities in the following cases.
\begin{enumerate}
\item $A=\Q$ and $n\geq5$;
\item $A=E$ and $n\geq\max\{3,|\mathcal O_E^\times|\}$;
\item $A=B$ and $n\geq4$.
\end{enumerate}
\end{cor}

Now we prove our main result.
\begin{cor}\label{cor:kodaira-degrees}
Let $\Gamma\subset H_{A,n}(\Q)$ be arithmetic and let $\ell$ be the weight
in Table~\ref{tab:cusp-weights} for $A$ and $n$.
If $\ell<k_{\mathrm{can}}$, then $X_\Gamma$ is of general type.
If $\ell=k_{\mathrm{can}}$, the cusp form in Theorem~\ref{thm:intro-cusp-forms}
gives a nonzero holomorphic differential of top degree on every smooth
projective model of $X_\Gamma$. In particular, $p_g(X_\Gamma)\geq1$
and $\kappa(X_\Gamma)\geq0$ in this case.
\end{cor}
\begin{proof}
Since $X_{\Gamma C_0}=X_\Gamma$, we may assume that $\Gamma$ contains $C_0$.
By Theorem~\ref{thm:intro-cusp-forms}, restriction to the component of
the identity gives a nonzero $F\in S_\ell(\Gamma)$.

Suppose first that $\ell<k_{\mathrm{can}}$. By Table~\ref{tab:cusp-weights},
$n$ satisfies the bound in Corollary~\ref{cor:canonical-toroidal}.
Let $X_\Gamma^{\mathrm{BB}}$ be the Baily--Borel compactification of $X_\Gamma$.
The corollary gives a projective toroidal compactification $\overline X_\Gamma^\Sigma$ with
canonical singularities and a natural morphism
$\overline X_\Gamma^\Sigma\to X_\Gamma^{\mathrm{BB}}$. We take the fans
$\Sigma$ regular for the full lattices of translations, with each cone
fixed pointwise by its stabilizer, as in Section~\ref{sec:boundary}.

For a rational boundary component $F_0$, let $U_{F_0}$ be the center of the
unipotent radical of its stabilizer. The quotient
$q:H_{A,n}\to H_{A,n}^{\mathrm{ad}}$ identifies $U_{F_0}$ with its image.
If $\gamma\in\Gamma$ satisfies $q(\gamma)\in q(U_{F_0})(\Q)$, then
$\gamma=(zI)u$ for some $u\in U_{F_0}(\Q)$. For $A=\Q,B$ one has $z=\pm1$.
For $A=E$, both $z$ and $z^{-1}$ are algebraic integers since $\Gamma$
preserves a lattice. Hence $z\in\mathcal O_E^\times$. In all cases
$u\in\Gamma$, so $q(\Gamma\cap U_{F_0}(\Q))=q(\Gamma)\cap q(U_{F_0})(\Q)$.
Thus every cusp is regular, and the toroidal charts are unramified along
the boundary divisors by \cite[Theorem~1 and Lemma~2]{ma2022boundary}.

Let $D_\Sigma$ be the reduced boundary and $\mathcal L$ the pullback of the
ample automorphic $\Q$-line bundle of weight one on $X_\Gamma^{\mathrm{BB}}$.
There are no branch divisors in the interior (see \cite[Section~II.1]{HulekSankaranSurvey},
\cite[(8.2)]{2026_Horinaga-Maeda_sing_unitary}, and Theorem~\ref{thm:DIII-compactification}).
The Hirzebruch proportionality \cite[Proposition~3.4(a)]{mumford1977hirzebruch} gives
\[
 K_{\overline X_\Gamma^\Sigma}\sim_\Q k_{\mathrm{can}}\mathcal L-D_\Sigma
 \sim_\Q (k_{\mathrm{can}}-\ell)\mathcal L+\operatorname{div}(F)-D_\Sigma.
\]
Here $\operatorname{div}(F)\sim_\Q\ell\mathcal L$ since the determinant character
for $A=E$ has finite order on $\Gamma$. Cuspidality at the regular cusps gives
$\operatorname{div}(F)\geq D_\Sigma$. Since $\mathcal L$ is big and
$\ell<k_{\mathrm{can}}$, the canonical divisor is big. The compactification
has canonical singularities, so $X_\Gamma$ is of general type.

If $\ell=k_{\mathrm{can}}$, the invariant form $F\omega$ descends to the smooth
unramified locus of $X_\Gamma$. Its $L^2$ norm is a constant multiple of the
Petersson norm of $F$, hence is finite. On a log resolution, the local coefficients
of $F\omega$ are square integrable and hence extend across the boundary divisor.
By birational invariance of holomorphic forms of top degree, this gives a nonzero
form on every smooth projective model. Its tensor powers give $\kappa(X_\Gamma)\geq0$.
\end{proof}

\subsection{Orthogonal modular varieties}

Let $(V,q)$ be a rational quadratic space of signature $(2,n)$,
where $n=2r\geq4$. We take $\DD_V$ to be a connected component of
$\{[z]\in\mathbb P(V_\C):q(z)=0,\ (z,\overline z)>0\}$,
where $(\ ,\ )$ is the bilinear form associated with $q$.
Let $\mathrm{GSO}(V)$ be the group of similitudes $g$ satisfying
$\det(g)=\mu(g)^{r+1}$, where $q(gv)=\mu(g)q(v)$.
We write $\mathrm{PGSO}(V)=\mathrm{GSO}(V)/\mathbb G_m$ for the quotient
over $\Q$ by the scalar matrices. Its real identity component is $\Aut(\DD_V)^\circ$.
Hilbert's theorem~90 gives $H^1(\Q,\mathbb G_m)=1$, so the quotient map induces
$\mathrm{GSO}(V)(\Q)/\Q^\times\simeq\mathrm{PGSO}(V)(\Q)$.
Fix an arithmetic subgroup $\GG\subset\mathrm{PGSO}(V)(\Q)$ preserving $\DD_V$.

\begin{thm}\label{thm:orthogonal}
For even $n\geq12$, the quotient $\GG\backslash\DD_V$ admits a projective toroidal
compactification with canonical singularities. It is of general type
if $n\geq22$, and has nonnegative Kodaira dimension if $n=16$ or $20$.
\end{thm}

The stabilizer at a point of $\DD_V$ is the image of $\SO(2)\times\SO(2r)$ in $\mathrm{PGSO}(V)(\R)$. 
Choose coordinates $e_0,e_1,\ldots,e_r$ on a maximal torus so that the weights on the
holomorphic tangent space are $e_0+e_i$ and $e_0-e_i$ for $1\leq i\leq r$.
The roots of the adjoint representation are
$\pm e_i\pm e_j$ for $0\leq i<j\leq r$.

\begin{prop}\label{prop:orthogonal-age}
If $n\geq12$, every nonidentity element
$s\in\mathrm{PGSO}(V)(\Q)$ of finite order fixing $p\in\DD_V$
satisfies $\age(s;T_p\DD_V)\geq1$.
\end{prop}
\begin{proof}
Conjugate $s$ into a maximal torus of the stabilizer.
Let $\alpha_i^+,\alpha_i^-\in[0,1)$ be the angles of the two eigenvalues on the tangent space corresponding to $e_0 \pm e_i$, and suppose that the age of $s$ is less than $1$.
The two weights $e_0+e_i$ and $e_0-e_i$ have sum $2e_0$, independently of $i$.
The sums $\alpha_i^++\alpha_i^-$ are therefore congruent modulo $\Z$.
They all lie in $[0,\age(s)]\subset[0,1)$, so each equals $\age(s)/r$. In particular, $0\leq\alpha_i^\pm\leq\age(s)/r<1/6$ since $r\geq6$.
Every compact root is a difference of two tangent weights, and the remaining roots are the tangent weights and their negatives.
Thus every eigenvalue of $\Ad(s)$ has real part greater than $1/2$. Since $\Ad(s)$ is defined over $\Q$ and has finite order, Lemma~\ref{lem:cyclotomic} gives $\Ad(s)=I$.
As $\mathrm{PGSO}(V)$ is adjoint, $s=1$, a contradiction.
\end{proof}

The following example shows the sharpness of the above bound.
\begin{ex}
Put
\[
 H=\begin{pmatrix}1&0\\0&3\end{pmatrix},\qquad
 M=\begin{pmatrix}3&-3\\1&3\end{pmatrix}.
\]
Direct calculation gives $M^tHM=12H$ and $M^6=-12^3I$.
Take $V$ with Gram matrix $H\oplus(-H)^{\oplus r}$ and set
$\gamma=M^{\oplus(r+1)}$. Then
$\gamma\in\mathrm{GSO}(V)(\Q)$ has multiplier $12$, and its image on $\mathrm{PGSO}(V)$ has order $6$.
The positive plane with Gram matrix $H$, oriented to give a point $x\in\DD_V$, is fixed by $[\gamma]$.
Choose $s=[\gamma]$ or its inverse so that its eigenvalues on $T_x\DD_V$
are $1$ with multiplicity $r$ and
$\exp(2\pi\sqrt{-1}/6)$ with multiplicity $r$.
Let $L$ be a lattice of $\Lie(G)(\bbQ)$ preserved by $s$.
For a torsion free normal arithmetic subgroup $\Gamma$ of the stabilizer of $L$, let $\Gamma'$ be the arithmetic subgroup generated by $\Gamma$ and $s$.
As in Example~\ref{ex:sharpness}, the stabilizer of $x$ in $\Gamma'$ is $\langle s\rangle$, since its stabilizer in $\Gamma$ is trivial.
The fixed space of each nonidentity power of $s$ in $T_x\DD_V$ has codimension $r\geq2$, so there are no quasireflections.
For $n=4,6,8,10$ we have $\age(s)=r/6=n/12<1$, so the Reid--Shepherd-Barron--Tai criterion gives a noncanonical singularity at the image of $x$.
\end{ex}

The examples use rational projective similitudes which do not lift to the orthogonal group. 
The bound $12$ therefore concerns a larger class of groups than the estimates for singularities of orthogonal modular varieties in
\cite{Gritsenko2007kodaira} and \cite[Theorem~2.5]{2026_Horinaga-Maeda_sing_unitary}.

For even $n\geq12$ put
$\GG_0=\GG\cap\operatorname{Im}(\SO(V)(\Q)\to\mathrm{PGSO}(V)(\Q))$.
Propositions~\ref{prop:orthogonal-age} and \ref{prop:central-descent}
show that the quotients by $\GG_0$ and $\GG$ have the same Kodaira dimension.

\begin{proof}[Proof of Theorem~\ref{thm:orthogonal}]
The assertion about canonical singularities follows from
Proposition~\ref{prop:orthogonal-age} and
Theorem~\ref{thm:boundary-transfer}.
The inverse image of $\GG_0$ in $\SO(V)(\Q)$ is arithmetic, so it
preserves a lattice $L\subset V$. After scaling the quadratic form
by a positive rational number, we may assume that $L$ is integral.
Thus $\GG_0\backslash\DD_V$ dominates
$\SO^+(L)\backslash\DD_V$, where $\SO^+(L)$ is the subgroup of $\SO(L)$ preserving $\calD_V$.
By \cite[Corollary~1.8]{ma2018kodaira}, this latter quotient is of
general type for even $n\geq22$, and its smooth projective models
have positive geometric genus when $4\mid n$ and $n\geq16$.
Hence the assertion holds.
\end{proof}

\appendix
\section{Automorphic representations without nontrivial self twist}
\subsection{Cusp forms coming from Hecke characters}
Let $E$ be an imaginary quadratic field with discriminant $-D$ and $\omega_{E/\bbQ}$ be the quadratic Hecke character corresponding to the quadratic extension $E/\bbQ$.
We denote the quadratic Dirichlet character with conductor $D$ by the same symbol $\omega_{E/\bbQ}$.

Let $F$ be a quadratic field with the nontrivial element $\sigma \in \Gal(F/\bbQ)$ and $\bbA_F$ be its adele ring.
Recall the cusp forms arising from a character $\chi$ of $F^\times \bs \bbA_F^\times$ with $\chi \neq \chi \circ \sigma$.
Let $\pi(\chi)$ be the automorphic induction in the sense of \cite{Langlands_Base_Change}.
Then $\pi(\chi) \isom \pi(\chi) \otimes (\omega_{F/\bbQ} \circ \det)$.
The cuspidal automorphic representation $\pi$ of $\GL_2(\bbA)$ satisfying $\pi \isom \pi \otimes (\eta \circ \det)$ with some nontrivial character $\eta$ can be written in this way.
By \cite[Proposition 11.1.12]{2021_computing_classical_MF}, if $\pi(\chi)$ corresponds to a holomorphic modular form of weight $k \geq 2$, then $F$ is imaginary.
Henceforth, we assume that $F$ is an imaginary quadratic field with discriminant $-d$.

Let $k>1$ be an odd positive integer.
If $\pi(\chi)$ corresponds to a modular form of weight $k$, then the infinity component of $\chi$ is of the form $z \mapsto (z/|z|)^{k-1}$ by replacing $\chi$ by $\chi \circ \sigma$ if necessary.
For a character $\mu$ of $F^\times \bs \bbA_F^\times$, let $\frakf_\mu$ be the conductor of $\mu$.
The corresponding modular form for $\pi(\chi)$ lies in $S_k(\Gamma_0(dN_{F/\bbQ}(\frakf_\chi)), \omega_{F/\bbQ}\chi')$.
Here, $\chi'$ is the Dirichlet character corresponding to $\chi|_{\bbA^\times}$ whose conductor divides $N_{F/\bbQ}(\frakf_\chi)$.

We consider the condition for the existence of a cuspidal automorphic representation of $\GL_2(\bbA)$ with no nontrivial self-twist and arising from a newform in $S_k(\Gamma_0(D), \omega_{E/\bbQ})$.
As discussed above, such a newform exists if and only if 
\[
\dim S_k(\Gamma_0(D), \omega_{E/\bbQ}) - \sum_{F,\chi} 1 \neq 0,
\]
where $F$ runs over all imaginary quadratic extensions with discriminant $-d$ and $\chi = \bigotimes_v' \chi_v$ runs over all characters of $F^\times \bs \bbA_F^\times$ such that $dN_{F/\bbQ}(\frakf_\chi) = D, \chi_\infty(z) = (z/|z|)^{k-1}$ and $\omega_{E/\bbQ} = \omega_{F/\bbQ} \chi'$.
The purpose of the appendix is to prove the following.

\begin{thm}\label{thm_existence_cusp_rep_no_twist}
Let $E$ have discriminant $-D$ and let $c_E$ be as in
Section~\ref{subsec:intro-kodaira}. For each odd integer $k\geq c_E+1$,
there is a newform in $S_k(\Gamma_0(D),\omega_{E/\Q})$ whose associated
cuspidal representation $\pi$ of $\GL_2(\bbA)$ has no nontrivial self twist.
        \end{thm}

\subsection{Proof of Theorem \ref{thm_existence_cusp_rep_no_twist}}

We first compute the number of $\chi$ that satisfy the assumption.
\begin{lem}\label{lem_appendix}
    Let $F$ be an imaginary quadratic field with discriminant $-d$ with $d \mid D$.
    If there exists a character $\chi$ of $F^\times \bs \bbA_F^\times$ such that $D= dN_{F/\bbQ}(\frakf_\chi), \chi_\infty(z) = (z/|z|)^{k-1}$ and $\omega_{E/\bbQ} = \omega_{F/\bbQ} \chi'$, the following conditions hold.
    \begin{itemize}
        \item $d$ and $D/d$ are relatively prime;
        \item $D/d \equiv 0,1 \bmod 4$;
        \item $(-d/p) = 1$ for any $p \mid (D/d)$, where $(\cdot/\cdot)$ is the Kronecker symbol.
            \end{itemize}
    Under these conditions, the number of such characters $\chi$ of $F^\times \bs \bbA_F^\times$ is equal to 
    \[
    2^{\omega(D/d)}h(F)f_k(d,D),\qquad
    f_k(d,D)
    =
    \begin{cases}
        1 & \text{if $d>4$;}\\
        1 & \text{if $d=4, (k-1)/2 \equiv (D/d-1)/4 \bmod 2$;}\\
        1 & \text{if $d=3, k \equiv 1 \bmod 6$;}\\
        0 & \text{otherwise.}
    \end{cases}
    \]
    where $\omega(D/d)$ is the number of distinct prime factors and $h(F)$ is the class number of $F$.
\end{lem}
\begin{proof}
    Let $c$ be the conductor of $\chi'$.
        Since $c$ divides $N_{F/\bbQ}(\frakf_\chi)$, we have 
    \[
    D \leq \mathrm{lcm}(d,c) \leq \mathrm{lcm}(d,N_{F/\bbQ}(\frakf_\chi)) = D/\gcd(d,D/d).
    \]
    Hence, $c = N_{F/\bbQ}(\frakf_\chi) = D/d$ and $(d,D/d) = 1$.
    Since $\chi'$ is a quadratic character and $\chi'(-1) = 1$, we have $D/d = c \equiv 0,1 \bmod 4$.
    Let $p\mid c$. Since $(d,c)=1$, the prime $p$ is unramified in $F$.
    If $p$ is inert in $F$, the conductor exponent $a$ of $\chi_p$ satisfies $\val_p(c) = 2a$, whereas $\val_p(c) \leq a$.
        This contradicts $p\mid c$.
    Thus every $p\mid c$ splits in $F$.
    
    For $p\mid c$, write $p\calO_F=\frakp\bar{\frakp}$ and let $s,t$
    be the conductor exponents at the two primes. The norm formula gives
    $s+t=\val_p(c)$. Restriction to the diagonal $\Z_p^\times$ gives
    $\val_p(c)\leq\max(s,t)$, so exactly one exponent is positive.
    Choose $\frakp$ to carry this exponent. Under $F_{\frakp}\simeq\Q_p$,
    the restriction of $\chi_{\frakp}$ to $\Z_p^\times$ is
    $\chi'_p|_{\Z_p^\times}$, and $\chi_{\bar{\frakp}}$ is unramified.
    Outside $c$ all finite components are unramified. This gives
    $2^{\omega(c)}$ choices for the restrictions to the finite local units.

    Together with $\chi_\infty$, each choice extends to a character
    if and only if it is trivial on $\calO_F^\times$. When an extension
    exists, the others are its twists by the characters of the ideal class
    group of $F$.
                            When $d>4$, we have $\calO_{F}^\times = \{ \pm 1\}$ and $\chi_\infty(-1) = +1$.
    By $D/d \equiv 0,1$, we have $\prod_v\chi_v(-1) = +1$.
    Hence, all such characters exist.
    When $d=4$, for $p \mid D/d$, define $x_p \in \bbF_p^\times$ by the image of $\sqrt{-1}$ in $\calO_{F} \rightarrow \calO_{F}/\frakp$.
    Then, $(x_p/p) = (-1)^{(p-1)/4}$ and $\prod_{v<\infty} \chi_v(\sqrt{-1}) \cdot \chi_\infty(\sqrt{-1}) = \prod_{p \mid (D/d)}(-1)^{(p-1)/4}\cdot (\sqrt{-1})^{(k-1)}$.
    Hence, $\chi$ is a character of $F^\times \bs \bbA_F^\times$ if and only if
    \[
    \frac{k-1}{2} \equiv\sum_{p \mid (D/d)} \frac{p-1}{4}
    \equiv \frac{D/d-1}{4}\bmod 2.
    \]
    When $d=3$, let $\zeta_3 = \exp(2\pi\sqrt{-1}/3)$.
        At each finite place the restriction of $\chi_v$ to $\calO_{F_v}^\times$ has order
    dividing two, so $\chi_v(\zeta_3)=1$.
    By the same method above, such $\chi$ is a character of $F^\times \bs \bbA_F^\times$ if and only if $k \equiv 1 \bmod 6$.
\end{proof}

\begin{proof}[Proof of Theorem \ref{thm_existence_cusp_rep_no_twist}]
    By \cite[Subsection 6.3]{2007_Stein}, we have
    \[
    \dim S_k(\Gamma_0(D), \omega_{E/\bbQ}) \geq \frac{k-1}{12} D\prod_{p\mid D}\left(1+\frac{1}{p}\right) - \frac{5 \cdot 2^{\omega(D)}}{6}.
    \]
            Hence, the statement holds for any $(k,D)$ satisfying 
    \[
    \frac{k-1}{12} D\prod_{p\mid D}\left(1+\frac{1}{p}\right) - \frac{5 \cdot 2^{\omega(D)}}{6}-\sum_{F}2^{\omega(D/d)}h(F)f_k(d,D) > 0.
    \]
    Here, $F$ runs over all imaginary quadratic fields with discriminant $-d$ such that $d$ satisfies the condition in Lemma \ref{lem_appendix}.
    It is well-known that $h(F) \leq \sqrt{d}(\log(d)+2)/\pi$ (cf.~\cite[p.49]{2000_Davenport}).
    Hence, we have
    \[
    \sum_F 2^{\omega(D/d)}h(F)f_k(d,D) \leq \sum_{d} \frac{2^{\omega(D)}\sqrt{d}(\log(d)+2)}{2^{\omega(d)}\pi} \leq \frac{\log(D)+2}{\pi}\sum_d\frac{2^{\omega(D)}\sqrt{d}}{2^{\omega(d)}}
    \]
    In the second and third terms, $d$ runs over all positive integers with $(d, D/d) = 1$ such that $d$ is the absolute value of a discriminant of an imaginary quadratic field.
    Also, $D\prod_{p \mid D}(1+p^{-1}) \geq \sqrt{D} \cdot 2^{\omega(D)}$.
        Thus it suffices to show that
    \[
    \frac{k-1}{12}>R(D),\qquad
    R(D)=\frac{5}{6\sqrt D}+\frac{2(\log D+2)}{\pi D^{1/3}}.
    \]
    Here we have used the bound
    \[
    \frac{\sum_d2^{\omega(D/d)}\sqrt{d}}{D\prod_{p\mid D}(1+p^{-1})}
    \leq \prod_{p\mid D}\frac{2+p^{\val_p(D)/2}}{p^{\val_p(D)-1}(p+1)}
    \leq 2D^{-1/3}.
    \]
    The product allows each full prime power dividing $D$ to occur in either $d$ or $D/d$.
    For the last inequality, multiply the factor at $p^a\parallel D$ by $p^{a/3}$.
    For odd $p$ we have $a=1$, and $x^{1/3}(2+\sqrt{x})/(x+1)$ decreases for $x\geq1$ and is less than one at $x=11$.
    If $2\mid D$, then $a=2$ or $3$, and the resulting factor at $2$ is at most $2^{5/3}/3$.
    The product of this bound and the factors at $3,5,7$ is less than $1.92$.
                                    The function $R(D)$ decreases for $D\geq3$. Since $R(3)<2$ and $R(10^6)<1/6$,
    the inequality holds for $k\geq25$ or $D\geq10^6$.

        For negative fundamental discriminants $-D$ with $D<10^6$ and $3\leq k\leq23$, we compared the dimension of $S_k(\Gamma_0(D),\omega_{E/\Q})$ from \cite[Section~6.3]{2007_Stein}
    with the dimension of the subspace given by Lemma~\ref{lem_appendix}.
        The dimensions agree exactly for $D\in\{3,4,7,8,11,15,20,23\}$ at weight $3$,
    for $D\in\{3,4,7\}$ at weight $5$, and for $D=3$ at weight $7$.
        \end{proof}

\bibliographystyle{alpha}
\bibliography{main}
\end{document}